\documentclass[12pt,reqno]{amsart}
\newcommand{\dist}{\mathop{\rm dist\,}}

\newcommand{\sgn}{\mathop{\rm sgn\,}}
\newcommand{\im}{\mathop{\rm Im\,}}
\newcommand{\re}{\mathop{\rm Re\,}}
\newcommand{\tcc}{the Cart\-wright class}
\newcommand{\tkc}{the Krein class}
\newcommand{\als}{Akhiezer--Levin set}
\newcommand{\mar}{{\mathcal M_E}}
\newcommand{\PI}{{PI}}

\newcommand{\pol}{{\mathcal P}}
\newcommand{\e}{{\varepsilon}}

\newtheorem{thm}{Theorem}[subsection]
\newtheorem{cor}[thm]{Corollary}
\newtheorem{lem}[thm]{Lemma}
\newtheorem{pro}[thm]{Proposition}

\newtheorem*{tha}{Theorem A}
\newtheorem*{thb}{Theorem B}
\newtheorem*{thc}{Theorem C}
\newtheorem*{thd}{Theorem D}
\newtheorem*{thee}{Theorem E}
\newtheorem*{thef}{Theorem F}

\theoremstyle{remark}
\newtheorem{rmk}[thm]{Remark}

\numberwithin{equation}{subsection}

\newcommand{\C}{{\mathbb C}}

\newcommand{\R}{{\mathbb R}}
\newcommand{\Z}{{\mathbb Z}}

\title[\tiny{Krein's entire functions}]
{Krein's entire functions and \\ the Bernstein approximation
problem}
\author[\tiny{Alexander Borichev and 
Mikhail Sodin}]{Alexander Borichev and 
Mikhail Sodin}
\thanks
{The second named author was supported by the Israel Science
Foundation of the Israel Academy
of Sciences and Humanities under Grant No. 93/97-1.}
\date{\today}

\begin{document}

\begin{abstract}
We extend two theorems of Krein concerning entire functions of
Cartwright class, and 
give applications for the Bernstein weighted
approximation problem.
\end{abstract}
\maketitle
\bigskip

\subsection{The Krein class and functions of bounded type.}

We start with two classical theorems of Krein concerning entire
functions. An entire function $f$ belongs to {\emph \tcc} if
$f$ has at most exponential type, that is
$$
\log|f(z)|=O(|z|),\qquad |z|\to\infty,
$$
and the logarithmic integral converges:
$$
\int_{-\infty}^{\infty}\frac{\log^+|f(x)|}{1+x^2}dx<\infty.
$$

\begin{tha} {\rm(Krein \cite{Kre1})} 
An entire function $f$ belongs to \tcc\
if and only
if the function $\log^+|f|$ has {\rm (}positive{\rm )} harmonic majorants in
both the upper and the lower half-planes.
\end{tha}

An entire function $f$ belongs to {\emph \tkc} 
if its zeros $\lambda_n$ are (simple and) real, and $1/f$ is  
represented as an absolutely convergent sum of simple fractions
\begin{equation}
\label{e1}
\frac1{f(z)}=\sum_{n}\frac 1{f'(\lambda_n)(z-\lambda_n)}\,, \qquad
\sum_{n}\frac 1{|f'(\lambda_n)|}<\infty\,.
\end{equation}

\begin{thb} {\rm(Krein \cite{Kre1})} 
The Krein class is contained in \tcc. 
\end{thb}

For the proofs see also \cite{Lev1}.
These two results have numerous applications in operator theory
and harmonic analysis (see, for example, \cite{Kre2}, \cite{Kre3}, 
\cite{deB2}, \cite[Chapter~IV]{GK}, \cite[Section~VI~F]{Koo}). 
Later, they were generalized in different directions (cf.
\cite[Section~26.4]{Lev1}, 
\cite[Section~VI.2]{GO}).

Let $E$ be a non-empty closed subset of the real line. In what
follows we assume that $E$ is regular for the Dirichlet problem
in $\C\setminus E$. A function $f$ analytic in $\C\setminus E$
is said to be \emph{of bounded type} if $\log^+|f|$ has a
harmonic majorant in $\C\setminus E$. It is well known that if
$f$ and $g$ are of bounded type, and $f/g$ is analytic in
$\C\setminus E$, then $f/g$ is also of bounded type there (cf.
\cite[Chapter~VII]{Nev}, \cite[Theorem~19, p.~181]{Par}). 
It is worth to mention that any function $\varphi$ lower
semicontinuous in the plane, which has a positive harmonic
majorant in $\C\setminus E$, satisfies the inequality
\begin{equation}
\int_E\varphi^+(x)\,\omega(i,dx,\C\setminus E)<\infty.
\label{12a}
\end{equation}

We would like to know when the assertions in two Krein's
theorems can be improved to guarantee that $f$ is of bounded
type in $\C\setminus E$. Note that every polynomial is of
bounded type in $\C\setminus E$. Indeed, our conditions on $E$
imply that $E$ has positive capacity, and therefore the identity
function omits in $\C\setminus E$ values from a set of positive
capacity. Hence, by the Frostman theorem, (see \cite[Chapter X,
Section~2.8]{Nev} for the case of the unit disc, and use the
uniformization argument in the general case), the identity
function is of bounded type in $\C\setminus E$, and the
statement for polynomials follows immediately.

For every regular set $E\subset \R$, we denote by $\mar (z)$ the
symmetric Martin function for $\C\setminus E$ with singularity
at infinity, that is, a positive harmonic function in
$\C\setminus E$ which vanishes on $E$ and satisfies the equality 
$\mar (\bar z)=\mar (z)$. A uniqueness theorem proved by
Benedicks \cite[Theorems~2 and 3]{BA} and Levin
\cite[Theorem~3.2]{LTF} asserts that $\mar $ exists and is
unique up to a positive multiplicative constant.  The function
$\mar$ extended by zero on $E$ is subharmonic in $\C$, and has
at most order one and mean type: $\mar(z)=O(|z|)$,
$|z|\to\infty$.
 
Our first result describes the sets $E$ such that every
Cartwright class function is of bounded type in $\C\setminus E$.
We say that a set $E\subset\R$ is {\emph {an \als }} if the
function $\mar$ is of mean type with respect to order $1$, that
is,
$$
\sigma_{\mar} \stackrel{\rm def}= \limsup_{|z|\to\infty} \frac{\mar(z)}{|z|}
>0\,.
$$ 
It is worth mentioning that in this case the limit 
$$
\sigma_{\mar} = \lim_{|y|\to\infty} \mar(iy)/|y|
$$ 
exists, and $\mar (z) \ge \sigma_{\mar} |\im z|$.
The function $\mar$ normalized by the condition $\sigma_{\mar}=1$ is
sometimes called the Phragm\'en--Lindel\"of function. 

The class of {\als}s was introduced in \cite{AL}. Let us present
two equivalent conditions. A set $E\subset\R$ is an
\als\ if and only if either of the following two properties holds:

\begin{enumerate}
\item (Koosis \cite[Section~VIII A.2]{Koo}) $\int_\R G(t,z)\,dt<\infty$,
where $G$ is the Green function for $\C\setminus E$, $z\in\C\setminus E$.
\item (Benedicks \cite[Theorem~4]{BA}) $\int_\R
\beta_E(t)/(1+|t|)\,dt<\infty$, where $\beta_E(t)$ 
is the harmonic measure $\omega(t,\partial
S_t,S_t\setminus E)$ of the boundary of the square
$S_t=\{z=x+iy:|x-t|<t/2,|y|<t/2\}$ with respect to the domain
$S_t\setminus E$ at the point $t$.
\end{enumerate}

Next we present three metric tests:

\begin{enumerate}
\item (Akhiezer--Levin \cite[Section~3.VII]{AL},
Kargaev \cite[Theorem~6(a)]{Kar}) If
$\int_{\R\setminus E}dx/(1+|x|)<\infty$,
then $E$ is an \als.
\item (Schaeffer \cite[Lemma 1]{Sch}) If $E$ is
relatively dense with respect to
Lebesgue measure $dm$ (that is, for some $a,b$, and for every
$x\in\R$, we have $m(E\cap[x,x+a])\ge b$), then $E$ is an \als.
\item (Kargaev \cite[Theorem~4]{Kar}) If $E$ is an \als, then
$\int_\R \dist(x,E)/(1+x^2)dx<\infty$.
\end{enumerate}

Given a positive symmetric harmonic function $h$ on $\C\setminus E$, set 
$$
C=\max\{c\ge 0: h-c\mar\text{\ is non-negative on\ }\C\setminus E\},
$$
and define the function $\PI_{E,h}$ (the Poisson integral of a
non-negative measure with support on $E$) which is non-negative,
symmetric, and harmonic on $\C\setminus E$, and satisfies the
equality
\begin{equation}
h=\PI_{E,h}+C\mar;
\label{61}
\end{equation}
clearly,
\begin{equation}
\text{there is no $\varepsilon>0$ such that\ } 
\PI_{E,h}\ge \varepsilon\mar \text{\ on\ } \C\setminus E. 
\label{161}
\end{equation}
\smallskip

We use a lemma which is possibly known. Since we were unable to
find the appropriate reference, we give its proof in Section~3.

\begin{lem} For every positive symmetric
harmonic function $h$ on $\C\setminus E$,
$$
\PI_{E,h}(iy)=o(\mar(iy)),\qquad |y|\to\infty.
$$
\label{pro*}
\end{lem}

Now, we present an extension of Theorem~A:

\begin{thm} If $E\subset\R$ is an \als, then every
function $f$ in \tcc\ is of bounded type in $\C\setminus E$.
Conversely, let $f$ be an entire function of non-zero exponential
type belonging to \tcc. If $f$ is of bounded type in $\C\setminus E$,
then $E$ is an \als.
\label{te1}
\end{thm}

\begin{proof} Let $E$ be an \als, and let $f$ be in \tcc, and of
exponential type $\sigma\ge 0$. First we suppose that $|f(x)|\le
1$, $x\in\R$. Applying the Phragm\'en--Lindel\"of principle to the 
function $\log|f|-\sigma_1\mar$ with $\sigma_1>\sigma\,
\sigma_{\mar}^{-1}$, in the upper and in the lower half-planes,
we conclude that $\log|f|-\sigma_1 \mar$ is non-positive
everywhere in $\C$, and therefore, $\sigma_1 \mar$ is a positive
harmonic majorant for $\log|f|$ in $\C\setminus E$.

In the general case, we use the Beurling--Malliavin multiplier
theorem \cite{BM}: there exists a function $g$ in \tcc\ with
$(1+|f(x)|)|g(x)|\le 1$, $x\in\R$. Applying the previous
argument, we obtain that $g$ and $fg$, and hence, $f$, are of
bounded type in $\C\setminus E$.

Now, let $f$ be an entire function of non-zero exponential type
belonging to \tcc. Suppose that $f$ is of bounded type in
$\C\setminus E$. 
Then the function
$\log|f(z)|$ has a positive harmonic majorant $h(z)$;
without loss of generality we may assume that h is symmetric,
$h(z)=h(\overline{z})$. By formula (\ref{61}), 
$h=\PI_{E,h}+C\mar$. 
Since the function $f$ has non-zero exponential type,
Lemma~\ref{pro*} implies that $\mar(iy)\ge c|y|$
for large $|y|$.
This implies that $E$ is an \als.
\end{proof}

\medskip
Our next result extends Theorem~B. Now we assume that
$E$ is the union of disjoint closed intervals $I_m=[a_m,b_m]$ with
$\dist(0,I_m)\to\infty$. Given an interval $I$ we denote
its length by $|I|$.

\begin{thm} If $f$ is a Krein class function with zeros
$\lambda_n$ on $E=\cup I_m$, and $|I_m|\ge c\dist(0,I_m)^{-M}$,
$c>0$, $M<\infty$, then $f$ is of bounded type in $\C\setminus E$.
\label{te2}
\end{thm}

\begin{proof} We are to prove that the function 
$\sum_{n} 1/[f'(\lambda_n)(z-\lambda_n)]$
is of bounded type in $\C\setminus E$.
Multiplying, if necessary, $f$ by a polynomial
with real zeros, we obtain
\begin{equation}
\label{e*}
\sum_n \frac{1+|\lambda_n|^{M}}{|f'(\lambda_n)|} <\infty\,.
\end{equation}
Without loss of generality, we may assume that 
$f'(\lambda_n)$ are real. Furthermore,
\begin{equation}
\sum_{n}\frac 1{f'(\lambda_n)(z-\lambda_n)}=\sum_{j=1}^2 g_j(z)= 
\sum_{j=1}^2\sum_{n}\frac {c_{n,j}}{z-\lambda_n},
\label{e2}
\end{equation}
where $c_{n,1}\ge 0$, $c_{n,2}\le 0$, and 
$$
\sum_{j=1}^2\sum_{n}(1+|\lambda_n|^{M})|c_{n,j}|<\infty.
$$
It suffices to verify that each $g_j$ in (\ref{e2}) is a
function of bounded type in $\C\setminus E$. Take one of them,
say $g_1$, and represent it as the sum of two functions,
$$
g_1(z)=g_-(z) + g_+(z)=
\sum_{\lambda_n\in E_-}\frac {c_{n,1}}{z-\lambda_n}+
\sum_{\lambda_n\in E_+}\frac {c_{n,1}}{z-\lambda_n},
$$
where $E_-=\bigcup[a_m,(a_m+b_m)/2)$,
$E_+=\bigcup[(a_m+b_m)/2,b_m]$, $E=E_+\cup E_-$.
Let us verify that $g_+$ is of
bounded type in $\C\setminus E$. Indeed, this function is
analytic in $\C\setminus E$ and satisfies the following
properties:
$$
\frac{\im g_+(z)}{\im z}<0,\qquad z\in\C\setminus\R,
$$
and, for $x\in \R\setminus E$,
$$
g_+(x)\ge -\sum_{\lambda_n\in E_+,\,\lambda_n>x}\, 
\frac{2c_{n,1}} {b_m-a_m}\ge -\frac{2}{c}
\sum_{n}c_{n,1}(1+|\lambda_n|^M)=\tau>-\infty\,.
$$
Thus, the image $g_+(\C\setminus E)$ omits the ray
$(-\infty,\tau)$. Applying the Frostman theorem, we get that
the function $g_+$ is of bounded type in $\C\setminus E$.
(Alternatively, we could just consider the functions
$\theta(z)=(z-1)/(z+1)$, $\theta_1(z)=(1+z)/(1-z)$, 
$\psi=\theta(\sqrt{g_+-\tau})$; $\psi$
has absolute values bounded by one in $\C\setminus E$. Hence,
$g_+=[\theta_1(\psi)]^2+\tau$ is of bounded type there.) The
same argument works for $g_-$, and we conclude that
$g_1=g_-+g_+$ is of bounded type in $\C\setminus E$, $g_2$, and
finally, $1/f$ are of bounded type in $\C\setminus E$.
\end{proof}

\begin{rmk} Using more information on \tkc\ function $f$ we can
further weaken our conditions on $|I_m|$: the assertion of
Theorem~\ref{te2} holds if for some $c>0$, $M<\infty$, and for every
zero $\lambda_n$ of $f$,
$$
|I(\lambda_n)|\ge \frac c{(1+|\lambda_n|^M)|f'(\lambda_n)|},
$$
where $I(\lambda_n)$ is the interval of $E$ containing the point
$\lambda_n$. It looks plausible that the assertion of
Theorem~\ref{te2} holds for any system of intervals of
``non-quasianalytically decaying lengths''. Namely, we may
conjecture that if $f$ is an entire function in \tkc, of
positive exponential type, with zeros $\lambda_n$, and if
$\varphi$ is a positive $\text{\rm Lip}\,1$ function on $\R$,
then $f$ is of bounded type in $\C\setminus E$, with
$$
E=\bigcup_n\left[\lambda_n-\frac 1{\varphi(\lambda_n)},
\lambda_n+\frac 1{\varphi(\lambda_n)}\right]
$$
if and only if $\int_\R \varphi(x)/(1+x^2)dx<\infty$. In some
special cases, when the zero set of $f$ is regularly
distributed, and $\varphi$ satisfies additional regularity
assumptions, this statement follows from Benedicks' results
\cite[Theorem~5]{BA}.
\end{rmk}

Combining Theorems~\ref{te1} and \ref{te2}, we get a sufficient
condition for $E$ to be an \als.

\begin{cor} If $f$ is a Krein class function 
of positive exponential type, with zeros $\lambda_n$, and 
if $M$ is a constant, then the set
$$
E=\bigcup_n\left[\lambda_n-\frac 1{(1+|\lambda_n|^M)|f'(\lambda_n)|}\,,
\lambda_n+\frac 1{(1+|\lambda_n|^M)|f'(\lambda_n)|}\right]
$$
is an \als.
\end{cor}
\bigskip

\subsection{The Bernstein problem on subsets of the real line.}

Fix a \emph{weight} $W$, that is a lower semicontinuous function
$W:\R\to[1,+\infty]$ such that $\lim_{|x|\to\infty}|x|^n/W(x)=0$
for $n\ge 0$. Consider the space $C(W)$ of functions $f$
continuous on $\R$ and such that
$\lim_{|x|\to\infty}|f(x)|/W(x)=0$. Put
$$
\|f\|_{C(W)}=\sup_{x\in\R}\frac{|f(x)|}{W(x)}.
$$
The Bernstein problem consists in finding out whether the set
$\pol$ of all polynomials is dense in $C(W)$. From now on, we
suppose that \emph{the weight $W$ is finite on a subset of $\R$
having a finite limit point}. In this case the polynomials are
simultaneously dense or not dense in every $C(W_r)$,
$W_r(x)=W(x)(1+|x|)^r$, $r\in \R$ (see
\cite[subsection~24]{Me}). Denote by $X_W$ the set of polynomials
$P$ such that $||P||_{C(W)}\le 1$. We define \emph{the Hall
majorant} $M_W$ as
$$
M_W(z)=\sup\{|P(z)|:P\in X_W\}.
$$
Since the function $\varphi$, $\varphi\equiv 1$, belongs to $X_W$, we have 
$M_W(z)\ge 1$, $z\in\C$. Furthermore, $M_W(x)\le W(x)$, $x\in\R$,
$\log M_W$ is lower semicontinuous in the plane.

General criteria of density of polynomials in weighted spaces
were obtained by Akhiezer--Bernstein, Pollard and Mergelyan at
the beginning of the 50-ies (cf. \cite{{Ak, Me, Koo}}). We shall use

\begin{thc} 
The polynomials are dense in $C(W)$ if and only if
one of three equivalent conditions holds: 
\begin{itemize}
\item[{\rm (1) }] {\rm{(Akhiezer--Bernstein)}}
$$
\sup_{P\in X_W}\, \int_{\R}\frac{\log^+|P(x)|}{1+x^2}\, dx =+\infty\,;
$$
\item[{\rm (2) }] {\rm{(Mergelyan)}}
$$
\int_\R\frac{\log M_W(x)}{1+x^2}\,dx=+\infty\,;
$$
\item[{\rm (3) }] {\rm{(Mergelyan)}} $M_W(z)=+\infty$
for some {\rm (}all\/{\rm )} $z\in\C\setminus\R$.
\end{itemize}
\end{thc}

Another criterion was proposed in \cite{deB1}.

\begin{thd} {\rm{(de~Branges)}}
The polynomials are not dense in $C(W)$ if and only
if there exists an entire function $F$ of zero exponential type,
$F\notin\pol$, with {\rm(}simple{\rm)} real zeros $\lambda_n$, such that
$$
\sum_n\frac{W(\lambda_n)}{|F'(\lambda_n)|}<+\infty.
$$
Such a function $F$ belongs to \tkc;
for such $F$ and for every polynomial $P$,
$$
\sum_n\frac{P(\lambda_n)}{F'(\lambda_n)}=0.
$$
\label{thd}
\end{thd}

The weight $W$ is assumed to be continuous in \cite{deB1};
however, the result holds for lower semicontinuous $W$ as well,
see \cite{SY}.

From now on, we suppose that $W(x)=\infty$ for $x\in\R\setminus
E$, where $E$ is a subset of $\R$ of the kind considered in the
previous section, $E=\cup I_m$, where $I_m$ are disjoint closed
intervals on $\R$, $\dist(0, I_m)\to\infty$ as $m\to\infty$.
Following Benedicks \cite{BP} and Koosis
\cite[Section~VIIIA]{Koo} we try to solve the Bernstein problem
for $C(W)$ in terms of $M_W{\bigm|}E$, 
replacing the form $(1+x^2)^{-1}dx$ by the harmonic measure
$\omega_E(dx)=\omega(i,dx,\C\setminus E)$.

\begin{thm} Suppose that $|I_m|\ge c(\dist(0,I_m))^{-M}$
for some $c>0$, $M<\infty$. The polynomials are dense
in $C(W)$ if and only if
\begin{equation}
\label{star}
\int_E\log M_W(x)\,\omega_E(dx)=+\infty.
\end{equation}
Furthermore, if the polynomials are not dense in $C(W)$, then
the function $\log M_W$ has a {\rm(}positive{\rm)} harmonic
majorant in $\C\setminus E$.
\label{te3}
\end{thm}

\begin{proof} To obtain the result in one direction, we prove that
\begin{equation}
\label{81}
\log|P(z)|\le \int_E\log^+|P(x)|\,
\omega(z,dx,\C\setminus E)\,, \qquad P\in\pol\,.
\end{equation} 
First, observe that our assumptions on $E$ imply that 
\begin{equation}
\label{beta}
\log |y| = o(\mar (iy))\,, \qquad y\to\infty\,.
\end{equation}
Indeed, take a sequence of points $x_k\in E$ tending to $\infty$
sufficiently rapidly (for example, $|x_{k+1}|>2|x_k|$ suffices),
and consider the entire function $B(z)=\prod_k (1-z/x_k)$. 
Then $B$ is in the Krein class, and $\log|B(iy)|/\log|y|
\to\infty$ as $y\to\infty$. By Theorem~\ref{te2}, $\log|B|$ has
a positive harmonic majorant $h$ in $\C\setminus E$. Using the
representation (\ref{61}) for $h$ together with Lemma~\ref{pro*}
we obtain (\ref{beta}). 

Then, applying a standard
Phragm\'en--Lindel\"of argument to the subharmonic functions
$$
\log |P(z)|-\int_E\log^+|P(x)|\,\omega(z,dx,\C\setminus E)
-c\mar(z)\,,\qquad c>0,
$$
in the domain $\C\setminus E$, we obtain (\ref{81}). 

By (\ref{81}),
$$
\log M_W(i)\le\int_E\log M_W(x)\,\omega_E(dx).
$$
If the last integral is finite, then $M_W(i)<\infty$, and by 
Theorem~C, the polynomials are not dense in $C(W)$.

Now we suppose that the polynomials are not dense in $C(W)$, and
as a consequence, are not dense in $C(W_r)$ where $r$ is to be
chosen later on. We apply Theorem~D, and get an entire function
$F$ in \tkc, of zero exponential type, $F\notin\pol$, with
zeros $\lambda_n\in E$, such that
\begin{equation}
\sum_n \frac{W_r(\lambda_n)}{|F'(\lambda_n)|}\le 1,
\label{63}
\end{equation}
and for every polynomial $P$, for every $z\in \C$,
$$
\sum_n \frac{P(z)-P(\lambda_n)}{(z-\lambda_n)F'(\lambda_n)}=0.
$$
Using relation (\ref{e1}) we get 
$$
\frac{P(z)}{F(z)}=
\sum_n \frac{P(\lambda_n)}{(z-\lambda_n)F'(\lambda_n)}\,.
$$
Theorem~\ref{te2} implies that the function $F$ is of
bounded type in $\C\setminus E$, and hence, by (\ref{12a}),
$$
0\le h(z)\stackrel{\rm def}=\int_E \log^+|F(t)|\,
\omega(z,dt,\C\setminus E)<\infty.
$$
Using this fact and relation (\ref{81}), we obtain
\begin{align}
\log |P(z)|\le &h(z)+\int_E \log^+\Bigl|
\sum_n \frac{P(\lambda_n)}{(t-\lambda_n)F'(\lambda_n)}\Bigr|\,
\omega(z,dt,\C\setminus E), \notag\\
\log  M_W(z)\le &h(z)+\int_E \log^+\sup_{P\in X_W} \Bigl|
\sum_n \frac{P(\lambda_n)}{(t-\lambda_n)F'(\lambda_n)}\Bigr|\,
\omega(z,dt,\C\setminus E),\notag
\end{align}
for $z\in\C\setminus E$, where $X_W$ is, as above, 
the set of polynomials $P$ such that $\|P\|_{C(W)}\le 1$. 
Therefore, to complete the proof of the theorem we need only to verify that
\begin{equation}
\int_E \log^+\sup_{P\in X_W} 
\Bigl|\sum_n \frac{P(\lambda_n)}{(t-\lambda_n)F'(\lambda_n)}\Bigr|\,
\omega_E(dt)<\infty.
\label{12b}
\end{equation}
Then, using Harnack's inequality, we conclude that $\log  M_W$
has a harmonic majorant in $\C\setminus E$, and by (\ref{12a})
we get that the condition (\ref{star}) does not hold.

Let us return to (\ref{12b}). By Jensen's inequality, for every
$d>0$ we have
\begin{gather}
\frac 1d \int_E 
\log^+\sup_{P\in X_W} 
\Bigl|\sum_n \frac{P(\lambda_n)}{(t-\lambda_n)F'(\lambda_n)}\Bigr|^d
\omega_E(dt)\notag\\
\le \frac 1d \log^+\int_E 
\sup_{P\in X_W} 
\Bigl|\sum_n \frac{P(\lambda_n)}{(t-\lambda_n)F'(\lambda_n)}\Bigr|^d
\omega_E(dt).\notag
\end{gather}
Furthermore, by (\ref{63}), for $0<d<1$ we get
\begin{align}
\sup_{P\in X_W} 
\Bigl|\sum_n\frac{P(\lambda_n)}{(t-\lambda_n)F'(\lambda_n)}\Bigr|^d 
&\le\sup_{P\in X_W} 
\sum_n\Bigl|\frac{P(\lambda_n)}{W_r(\lambda_n)}\Bigr|^d
\Bigl|\frac{W_r(\lambda_n)}{(t-\lambda_n)F'(\lambda_n)}\Bigr|^d
\notag\\&\le\sum_n \frac {1}{(1+|\lambda_n|)^{rd}|t-\lambda_n|^d}.
\notag
\end{align}
Therefore, to prove (\ref{12b}) it remains to check that for some $0<d<1$,
\begin{equation}
\sum_n\int_E \frac{1}{(1+|\lambda_n|)^{rd}|t-\lambda_n|^d}\,
\omega_E(dt)
<\infty.
\label{64}
\end{equation}

Since $\lambda_n$ are the zeros of an entire function of zero
exponential type, $1+|\lambda_n|\ge cn$ for some $c>0$. 
Thus, inequality (\ref{64}) follows from the estimate
\begin{equation}
\int_E \frac{1}{|t-\lambda|^d}\,\omega_E(dt)\le
C(1+|\lambda|)^{rd-2}
\label{65}
\end{equation}
with $C$ independent of $\lambda\in E$, and with some $0<d<1$.

Our conditions on $E$ imply that for some $c>0$, $M<\infty$,
and every $\lambda\in E$, there exists $\delta$,
$c(1+|\lambda|)^{-M}\le \delta\le 10 c(1+|\lambda|)^{-M}$\
such that for $I=[\lambda-\delta,\lambda+\delta]$ we have 
$$
E\cap I=J_1\cup J_2
$$
where intervals $J_k=[a_k,b_k]$ satisfy the condition
$|J_k|\ge c(1+|\lambda|)^{-M}$, $k=1,2$. 
The following elementary estimate of harmonic measure,
$$
\omega_E(dt)=
\omega(i,dt,\C\setminus E)\le \omega(i,dt,\C\setminus I_k)\le
\frac{C\,dt}{\sqrt{(b_k-t)(t-a_k)}},\quad t\in J_k,
$$
shows that for $d=1/4$,
\begin{multline}
\int_{J_k} \frac{1}{|t-\lambda|^d}\,\omega_E(dt)\\
\le\int_{a_k}^{b_k} \frac{C\,dt}
{(b_k-t)^{1/2}(t-a_k)^{1/2}|t-\lambda|^{1/4}}\le C(1+|\lambda|)^{M/4},
\quad k=1,2.\notag
\end{multline}
Furthermore,
$$
\int_{E\setminus I} \frac{1}{|t-\lambda|^{1/4}}\,\omega_E(dt)\le
\sup_{t\in E\setminus I}\frac{1}{|t-\lambda|^{1/4}}\le C(1+|\lambda|)^{M/4}.
$$
Thus, estimate (\ref{65}) is true for $r\ge M+8$, $d=1/4$.
\end{proof}

\begin{rmk} The same proof shows that the polynomials are not
dense in $C(W)$ as soon as (\ref{star}) fails and there exists a
sequence $m_k\to\infty$ such that $|I_{m_k}|\ge
c(\dist(0,I_{m_k}))^{-M}$ for some $c>0$, $M<\infty$.
\end{rmk}

\begin{rmk} In the general setting, if the polynomials are not
dense in $C(W)$, then every function $f$ from the closure of the
polynomials ${\rm Clos}_{C(W)}\pol$ has analytic continuation
into the whole complex plane, and $|f(z)|\le M_W(z)||f||_{C(W)}$,
$z\in \C$. Therefore, in the assumptions of Theorem~\ref{te3},
every function from ${\rm Clos}_{C(W)}\pol$ is of bounded type
in $\C \setminus E$.
\end{rmk}

\begin{rmk} As in Theorem~C, condition (\ref{star}) is equivalent to 
$$
\sup_{P\in X_W}\int_E\log^+|P(x)|\,\omega_E(dx)=+\infty\,.
$$
\end{rmk}
\bigskip

The following examples demonstrate that the assertions of
Theorem~\ref{te3} are not valid if the condition $|I_m|\ge
c(\dist(0,I_m))^{-M}$ is not fulfilled.

\begin{pro} {\rm(a)} There exist a weight $W$ and a subset $E$ of $\R$ 
such that 
$W(x)=\infty$ for $x\in\R\setminus E$, the polynomials are dense
in $C(W)$ and
$$
\int_E\log M_W(x)\,\omega_E(dx)<\infty.
$$
{\rm(b)} There exist 
a weight $W$ and a subset $E$ of $\R$ such that 
$W(x)=\infty$ for $x\in\R\setminus E$, the polynomials are not dense
in $C(W)$, and
$$
\int_E\log M_W(x)\,\omega_E(dx)=+\infty.
$$
\end{pro}

\begin{proof} (a) Consider a set of disjoint intervals $I_n$, such that 
$|I_n|\le 1$, $\exp n\in I_n$, $n\ge 1$, and 
$\omega(i,I_n,\C\setminus(I_1\cup I_n))<n^{-2}\exp(-n)$
for $n>1$. We define a weight $W$ as follows: 
$W{\bigm|}I_n\equiv\exp\exp n$, $n\ge 1$, and $W(x)=+\infty$ elsewhere.
By Theorem~D, the polynomials are not dense in $C(W)$. 
Indeed, no entire function $F$ of zero exponential type, 
with real zeros $\lambda_n\to\infty$, satisfies the condition
$$
|F'(\lambda_n)|\ge c\exp\lambda_n.
$$
Finally, since $M_W(x)\le W(x)$, $x\in\R$, we get
\begin{gather}
\int_E\log M_W(x)\,\omega_E(dx)
\le\!\int_{I_1}\log W(x)\,\omega_E(dx)+\sum_{n>1}
\int_{I_n}\log W(x)\,\omega_E(dx)\notag\\
\le C+\sum_{n>1}
\omega(i,I_n,\C\setminus(I_1\cup I_n))\sup_{I_n}\log W
\le C+\sum_{n>1}
n^{-2}e^{-n}e^n
<\infty.\notag
\end{gather}

\noindent
(b) We use an auxiliary estimate of harmonic measure.

\begin{lem} Let $I_n$ be intervals of length $1$ such that
$\dist(0,I_n)=(1+o(1))\exp n$, $n\to+\infty$, $E=\cup I_n$. Then
for some $\e>0$ and $C>0$, independent of $E$,
$$ 
\omega(i,I_n,\C\setminus E)\ge C\cdot e^{(\e-1)n}. 
$$
\label{le6}
\end{lem}

We postpone the proof of this lemma till the next section.

Fix $\rho$ with $1-\e<2\rho<1$, and consider the canonical product
$$
B_\rho(z)=\prod_{n=1}^\infty\Bigl(1-\frac z{n^{1/\rho}}\Bigr),
$$
and the entire function $F_\rho$ of zero exponential type,
$F_\rho(z)=B_\rho(z^2)$. Denote by $\Lambda$ the zero set of
$F_\rho$, $\Lambda=\bigl\{\lambda_{\pm n}=\pm
n^{1/(2\rho)},\,n\ge 1\bigr\}$. It follows from the results of
G.~H.~Hardy in \cite{Ha} that the expression
$$
\frac{B_\rho(z)}{z^{-1/2}\sin(\pi z^\rho)
\exp((\pi\cot\pi\rho)z^\rho)}
$$
tends to a finite non-zero limit for $|z|\to\infty$ with $|\arg
z|\le\pi/2$. Put $f(z)=\sin(\pi z^\rho)$, $g=B_\rho/f$. 
Since $B^\prime_\rho(\lambda)=g(\lambda)f'(\lambda)$,
$\lambda=n^{1/\rho}$, $n\ge 1$, we get for some $c>0$:
\begin{gather}
|B^\prime_\rho(\lambda))|=(c+o(1))\lambda^{\rho-3/2}
\exp((\pi\cot\pi\rho)\lambda^\rho),
\qquad\lambda=n^{1/\rho},\,n\to\infty,\notag\\
|F^\prime_\rho(\lambda)|=(2c+o(1))|\lambda|^{2\rho-2}
\exp((\pi\cot\pi\rho)|\lambda|^{2\rho}),
\qquad\lambda\in\Lambda,\,|\lambda|\to\infty,\notag
\end{gather}
whence
\begin{multline}
\log|F^\prime_\rho(\pm\exp t)|=\log(2c)+o(1)+
(2\rho-2)t+(\pi\cot\pi\rho)\exp(2\rho t),\\
\qquad\exp t\in\Lambda,\,t\to\infty.\notag
\end{multline}

We define an even weight $W$ by the formula 
\begin{equation}
\log W(\pm\exp t)=(\pi\cot\pi\rho)\exp(2\rho t).
\label{79}
\end{equation}
Since
$$
\sum_{\lambda\in\Lambda}\frac{W(\lambda)}
{|\lambda|^2|F^\prime_\rho(\lambda)|} <\infty,
$$
Theorem~D implies that the polynomials are not dense in
$C(W_r)$. The even log-convex function $W$ is increasing on the
positive half-line. Let us verify that
\begin{equation}
\log W(x)=(1+o(1))\log M_W(x), \qquad x\to+\infty.
\label{19}
\end{equation}
First, without loss of generality, we assume that $W$ is
$C^2$-smooth. By convexity of $\varphi=\log W\circ\exp$, for every
$r,s$ we have $\varphi(s)+\varphi'(s)(r-s)\le\varphi(r)$. Now we fix
$t=\log x$. If $\varphi'(r)=n\le\varphi'(t)<n+1$, then $\varphi(t)\ge
n(t-r)$. Furthermore, for some $\xi\in (r,t)$ we have
$\varphi(t)-\varphi(r)=\varphi'(\xi)(t-r)$. Since $\varphi'(\xi)\le
\varphi'(t)<n+1=\varphi'(r)+1$, we get $\varphi(t)-\varphi(r)-
\varphi'(r)(t-r)\le t-r$, and as a result,
$\varphi(r)+\varphi'(r)(t-r)\ge (1-1/n)\varphi(t)$. Therefore,
$x^n\exp(\varphi(r)-rn)\ge W(x)^{1-1/n}$; for every $y>0$,
$r\in\R$, we have $W(y)\ge y^n\exp(\varphi(r)-rn)$, thus the
estimate (\ref{19}) is proved.

Relations (\ref{79}) and (\ref{19}) imply that
\begin{equation}
\log M_W(x)=(1+o(1))(\pi\cot\pi\rho)|x|^{2\rho},\qquad |x|\to\infty.
\label{94}
\end{equation}

Choose a sequence of intervals $I_n$ with $|I_n|=1$,
$\dist(0,I_n)=(1+o(1))\exp n$, $n\to+\infty$, and set $E=\cup
I_n$. Apply Lemma~\ref{le6}, and add to $E$ a union $E'$ of
small intervals such that $\Lambda\subset E'$, and
$$ 
\omega\bigl(i,I_n,\C\setminus (E\cup E')\bigr)\ge C_1\cdot e^{(\e-1)n}. 
$$
Then by (\ref{94}),
\begin{gather*}
\int_{E\cup E'}\log M_W(x)\,\omega\bigl(i,dx,\C\setminus (E\cup E')\bigr)
\\ 
\\
\ge \sum_n
\omega\bigl(i,I_n,\C\setminus (E\cup E')\bigr)
\inf_{I_n}\log M_W
\ge \sum_n C_2\cdot e^{(\e-1)n}e^{2\rho n}=+\infty.
\end{gather*}
\end{proof} 

\medskip
As corollaries to Theorem \ref{te3} we obtain some results of
Levin and Bene\-dicks. First, suppose that $W{\bigm|}E=
\widetilde{W}{\bigm|}E\rule{0pt}{16pt}$ for an even log-convex
function $\widetilde{W}$ increasing on the positive half-line,
$W{\bigm|}\R{\setminus} E\equiv+\infty$. Applying
Theorem~\ref{te3} and using relation (\ref{19}) for
$\widetilde{W}\rule{0pt}{16pt}$, we arrive at a statement, essentially
equivalent to that proved by Levin in
\cite[Theorem~3.23]{Lev2}:

\begin{thee} If $W$ is a weight as above, and a set $E$
satisfies the conditions of Theorem~{\rm\ref{te3}}, then 
the polynomials are dense in $C(W)$ if and only if
$$
\int_E\log W(x)\,\omega_E(dx)=+\infty.
$$
\end{thee}

Benedicks investigated in \cite{BP} (see also the discussion in  
\cite[Section~VIII A.4]{Koo})
the weighted polynomial approximation on the sets 
\begin{equation}
E=\cup_{n\in\Z}\bigl[|n|^p\sgn n-\delta,|n|^p\sgn n+\delta\bigr]
\label{69}
\end{equation} 
for $p>1$, $\delta<1/2$. He announced the following result:

\begin{thef} Suppose that $E$ is a set of the form
{\rm(\ref{69})}, and $W$ is a weight such that $W(x)=+\infty$ for
$x\in\R\setminus E$. The polynomials are dense in $C(W)$ if and
only if
$$
\sup_{P\in X_W}\int_E\frac{\log |P(x)|}{1+|x|^{1+1/p}}\,dx=+\infty\,.
$$
\label{thef}
\end{thef}

In \cite{BP} Benedicks gave a proof of the ``only if'' 
part of this theorem based on his estimate of the harmonic
measure $\omega_E(dx)=\omega(i,dx,\C\setminus E)$ for sets $E$
of the form (\ref{69}):

\begin{lem} If $E$ has the form {\rm(\ref{69})}, then 
\begin{multline}
\frac {c}{1+|x|^{1+1/p}}
\frac {1}{\sqrt{\delta^2-(x-|n|^p\sgn n)^2}} \\
\label{ben}
\le  
\frac{\omega_E(dx)}{dx} 
\le
\frac {C}{1+|x|^{1+1/p}}
\frac {1}{\sqrt{\delta^2-(x-|n|^p\sgn n)^2}}\,, 
\end{multline}
for $x\in \bigl[|n|^p\sgn n-\delta,|n|^p\sgn n+\delta\bigr]$,
where $c$ and $C$ are positive constants that do not depend on
$x$ and $n$.
\label{le5}
\end{lem}

A more accessible reference for the upper bound in (\ref{ben})
is \cite[Section~VIII~A.4]{Koo} where the reader may find a
sketch of the proof. A proof for the lower bound is given in the
next section.

Our Theorem~\ref{te3} together with the lower bound in
(\ref{ben}) immediately yields the ``if'' part of Theorem~F:

\begin{cor} Suppose that $W$ is a weight and $E$ is 
a set of the form {\rm(\ref{69})} such that 
$W(x)=\infty$ for $x\in\R\setminus E$.
If the polynomials are not dense in $C(W)$, then
$$
\int_E\frac{\log M_W(x)}{1+|x|^{1+1/p}}\,dx<\infty.
$$
\end{cor}
\bigskip

\subsection{Harmonic estimation in slit domains.}

Here, we prove Lemmas~\ref{pro*} and \ref{le6}, and the lower estimate
in Lemma~\ref{le5}.

\begin{proof}[Proof of Lemma~{\rm\ref{pro*}}]
Suppose that for a sequence $y_k\to+\infty$ we have
$$
\infty>\PI_{E,h}(iy_k)\ge\mar(iy_k).
$$
By Harnack's inequality, for some positive $c_1,c_2$ independent
of $k$,
\begin{equation}
\left.
\begin{aligned}
\PI_{E,h}(z)&\ge c_1\PI_{E,h}(iy_k)\\
\mar(z)&\le c_2\mar(iy_k)
\end{aligned}
\right\},\qquad |z-iy_k|<y_k/2.
\label{100}
\end{equation}
For some positive $c_3$, consider the function 
$u=c_3\mar-\PI_{E,h}$ harmonic on $\C\setminus E$. Then
\begin{equation}
u(z)\le -(c_1-c_2c_3)\mar(iy_k),\qquad |z-iy_k|<y_k/2.
\label{111}
\end{equation}

Next, we use the following fact (cf. Lemma~1 of \cite{Sch},
Lemma~6 of \cite{BA}):
\begin{equation}
\text{the function}\quad y\to\mar(x+iy)\quad \text{is increasing for} 
\quad y\ge 0. 
\label{53}
\end{equation}
For the sake of completeness, we give an argument from \cite{LTF}.
Since $\mar$ is positive, harmonic and symmetric
in $\C\setminus E$, is subharmonic in the plane, 
and has at most order one and mean type there,
the subharmonic version of the Hadamard representation 
(see \cite[Section~4.2]{HK})
implies that for a finite positive measure $\mu$ on $\R$,
\begin{align}
\mar(z)=&\int_{|t|\ge 1}\Bigl(\log\Bigl|1-\frac zt\Bigr|+\frac{\re z}{t}\Bigr)
\,d\mu(t)\notag\\+&\int_{|t|\le 1}\log|t-z|\,d\mu(t)+c_1+c_2\re z,\qquad z\in 
\C\setminus E.\notag
\end{align}
Now, the property (\ref{53}) follows immediately.

Using (\ref{53}) and the second estimate in (\ref{100}), we get
\begin{equation}
u(z)\le c_3\mar(z)\le c_2c_3\mar(iy_k),\quad |\re z|=y_k/2,\, |\im z|<y_k.
\label{112}
\end{equation}

Denote by $H$ the union of two horizontal sides of the domain
$S=\{z\in\C:|\re z|<y_k/2, |\im z|<y_k\}$, and by 
$V$ the union of its two vertical sides. An estimate of 
harmonic measure in $S\setminus E$ gives us that for some positive $C$
independent of $k$ and $E$,
\begin{equation}
\frac{\omega(z,H,S\setminus E)}{\omega(z,V,S\setminus E)}\ge \frac 1C,
\qquad |z|< y_k/5.
\label{85}
\end{equation}
To obtain this estimate we use an easy generalization of a lemma of Benedicks 
\cite[Lemma~7]{BA} (see also \cite[p.436]{Koo}).
This generalization claims that for every square 
$S_{x,t}=\{z\in\C:|\re z-x|< t, |\im z|< t\}$ with 
horizontal sides $H_{x,t}$ and vertical sides $V_{x,t}$,
the following inequality holds:
\begin{equation}
\omega(x+iy,H_{x,t},S_{x,t}\setminus E)\ge
\omega(x+iy,V_{x,t},S_{x,t}\setminus E),\, x+iy\in S_{x,t}.
\label{84}
\end{equation}

To verify (\ref{84}) we note first that by symmetry,
on the diagonals of the square $S_{x,t}$ we have
$\omega(\cdot,H_{x,t},S_{x,t})=\omega(\cdot,V_{x,t},S_{x,t})$.
Therefore, applying the maximum principle to the difference of
these functions, we get 
\begin{gather*}
\omega(r,H_{x,t},S_{x,t})\le \omega(r,V_{x,t},S_{x,t}),\qquad r\in(x-t,x+t),\\
\omega(x+iy,H_{x,t},S_{x,t})\ge \omega(x+iy,V_{x,t},S_{x,t}),
\qquad x+iy\in S_{x,t}.
\end{gather*}
Finally,
\begin{multline*}
\omega(x+iy,H_{x,t},S_{x,t}\setminus E)\\
=\omega(x+iy,H_{x,t},S_{x,t})-
\int_E \omega(r,H_{x,t},S_{x,t})\omega(x+iy,dr,S_{x,t}\setminus E)\\
\ge\omega(x+iy,V_{x,t},S_{x,t})-
\int_E \omega(r,V_{x,t},S_{x,t})\omega(x+iy,dr,S_{x,t}\setminus E)\\
=\omega(x+iy,V_{x,t},S_{x,t}\setminus E),\qquad x+iy\in S_{x,t}.
\end{multline*}

To deduce (\ref{85}) note that the function
$\omega(z,H,S\setminus E)$ is positive and continuous in
$S\setminus E$. Therefore, for $A=\{x\pm iy_k/5:|x|\le 2y_k/5\}$
we have $\min_{z\in A} \omega(z,H,S\setminus E)>0$. Hence, for sufficiently
big $C$,
$$
\varphi(z)\stackrel{\rm def}=C \omega(z,H,S\setminus E)-
\omega(z,V,S\setminus E)\ge 1, \qquad z\in A.
$$
For every $z=x+iy$ with $|z|< y_k/5$ consider the square
$S_{x,t}$ with $t=y_k/5$, and note that $H_{x,t}\subset A$,
$S_{x,t}\subset S$. Therefore, $\varphi{\bigm|}H_{x,t}\ge 1$,
$\varphi{\bigm|}E\equiv 0$, $\varphi{\bigm|}V_{x,t}\ge -1$, and the estimate
(\ref{84}) implies that $\varphi(x+iy)\ge 0$. Thus, the property
(\ref{85}) is proved.

Now, if $c_3$ is sufficiently small, then applying the theorem
on two constants to the symmetric harmonic function $u$ in the
domain $S\setminus E$, and using (\ref{111}), (\ref{112}),
(\ref{85}), and the property
$$
\limsup_{z\to w} u(z)=-\liminf_{z\to w} \PI_{E,h}(z)\le 0, \qquad w\in E, 
$$
we obtain 
$$
u(z)\le 0,\qquad |z|< y_k/5.
$$
Thus, $c_3\mar(z)-\PI_{E,h}(z)=u(z)\le 0$, $z\in\C\setminus E$. 
Applying the property (\ref{161}) we come to a contradiction.
\end{proof}

\begin{proof}[Proof of the lower estimate in Lemma~{\rm\ref{le5}}]

Step A. For $t\ge 0$ denote by $K_t$ the square $\{z\in\C:|\re
z-t|\le t/2, |\im z|\le t/2\}$. In what follows we use a
function $u(z)=\log|z+\sqrt{z^2-1}|$; it is positive and
harmonic in $\C\setminus[-1,1]$, vanishes on $[-1,1]$, and
$u(z)=\log|z|+O(1)$, $|z|\to\infty$.

The function $v(z)=u(Ct^{1-1/p}\sin\pi z^{1/p})$ vanishes on
a closed set $F$, $K_t\cap E\subset F$,
for some $C$ independent of $t$. Furthermore,
$v$ is non-negative on $K_t$, and harmonic on $K_t\setminus F$.
We estimate the function $z\to |t^{1-1/p}\sin\pi z^{1/p}|$
as follows:
\begin{gather}
|t^{1-1/p}\sin\pi z^{1/p}|\le \exp (Ct^{1/p}),\qquad z\in K_t,\notag\\
|t^{1-1/p}\sin\pi t^{1/p}|\ge Cn^{p-1},\qquad t=(n+1/2)^p,\,n\ge 0.\notag
\end{gather} 
The asymptotical relation $u(z)=\log|z|+O(1)$, $|z|\to\infty$,
implies now that $v(z)\le ct^{1/p}$, $z\in K_t$, and
$v((n+1/2)^p)\ge c\log n$, $n\ge 0$.  Next, we choose
$t=(n+1/2)^p$, and apply the theorem on two constants to the
function $v(z)$ in $K_t\setminus F$:
$$
c\log n\le v((n+1/2)^p)
\le \omega(t,\partial K_t,K_t\setminus F)\,\sup_{\partial K_t}v\le
c\,n\,\omega(t,\partial K_t,K_t\setminus F)\,.
$$
Hence, for $t=(n+1/2)^p$, $n>1$,
$$ 
\omega(t,\partial K_t,K_t\setminus E)\ge 
\omega(t,\partial K_t,K_t\setminus F)\ge c\,\frac {\log n}n\,.
$$

Let $H_t$ be the union of two horizontal sides of $K_t$. The
lemma of Benedicks mentioned in the proof of Lemma~\ref{pro*}
claims that
$$
\omega(t,H_t,K_t\setminus E)\ge 
\frac 12\omega(t,\partial K_t,K_t\setminus E).
$$
Therefore, for $n>1$,
\begin{equation} 
\omega((n+1/2)^p,H_{(n+1/2)^p},K_{(n+1/2)^p}\setminus E)\ge 
c\frac {\log n}{n}.
\label{14}
\end{equation}

Step B. The Green function $G(z,i)$ for $\C\setminus E$ is
positive, bounded and harmonic on $\{z:|\im z|>2\}$. Therefore,
applying the Poisson formula in the half-planes $\{z:\pm\im
z>2\}$ we get
\begin{multline}
G(z,i)\ge \frac1\pi\int_{-\infty}^\infty
G(x\pm 2i,i)\frac{|\im z|-2}{(|\im z|-2)^2+x^2}\,dx\\
\ge\frac{c}{|\im z|},\qquad 3|\im z|\ge |\re z|\ge 10.
\label{15}
\end{multline}
The inequalities (\ref{14}) and (\ref{15}) imply that for $n>1$,
$$
G((n+1/2)^p,i)\ge 
\omega((n+1/2)^p,H_{(n+\frac{1}{2})^p},K_{(n+\frac{1}{2})^p}\setminus
E)\inf_{H_{(n+\frac{1}{2})^p}}G 
\ge c\,\frac {\log n}{n^{p+1}}. 
$$
Set $I_n=[n^p-\delta,n^p+\delta]$.
Since $G(z,i)$ is positive and harmonic on $\{z:|z-n^p|\le n^{p-1}\}
\setminus I_n$, by Harnack's inequality we get
$$
G(z,i)\ge c\frac {\log n}{n^{p+1}},\qquad |z-n^p|=n^{p-1}.
$$

Step C. Finally, we consider an auxiliary function $w(z)=u((z-n^p)/\delta)$. It
is harmonic on $\C\setminus I_n$, vanishes on $I_n$, and
$w(z)\le c\log n$ for $|z-n^p|=n^{p-1}$. Therefore, for $n>1$,
$$
G(z,i)\ge \frac {c}{n^{p+1}}w(z), \qquad |z-n^p|\le n^{p-1},
$$
and
\begin{multline}
\frac{\omega(i,dx,\C\setminus E)}{dx}=\frac 1{\pi}
\frac{\partial G(x+iy,i)}{\partial y}\Bigm|_{y=0}\\
\ge\frac {c}{1+|x|^{1+1/p}}
\frac {1}{\sqrt{\delta^2-(x-|n|^p\sgn n)^2}}\,dx,\qquad
x\in I_n.\notag
\end{multline}
The estimate for $n\le 1$ is obtained in an analogous way.
\end{proof}

\begin{proof}[Proof of Lemma~{\rm\ref{le6}}]
As in the previous proof, let us consider 
the Green function $G(z,i)$ for $\C\setminus E$. It is positive, 
harmonic and bounded in $\C\setminus (E\cup\{z:|z-i|<1\})$.
Denote 
$$
h_n=\int_0^{\exp n}G(x,i)\,dx,\qquad n\ge 0.
$$
Applying the Poisson formula in the lower half-plane, we get 
$$
G(-ie^n,i)\ge \frac1\pi\int_0^{\exp n}\frac{e^n}{x^2+e^{2n}}G(x,i)\,dx
\ge c\cdot e^{-n} h_n,\qquad n\ge 0.
$$
By Harnack's inequality, at least on one half of the length of
the interval $[e^n,e^{n+1}]$, we have $G(x,i)\ge c\cdot
G(-ie^n,i)\ge c\cdot e^{-n} h_n$. Therefore, for some $c>0$, we
have $h_{n+1}>(1+c)h_n$, $n\ge 0$. Consequently,
$h_{n+1}>c(1+c)^n$, $n\ge 0$. Once again, by the Poisson
formula, for some $\e>0$,
$$
G(-ix,i)>c\cdot x^{\e-1},\qquad x\ge 1.
$$
For every $n$ we denote by $c_n$ the center of $I_n$. Arguing as
in step C of the previous proof, we compare $G(z,i)$ with
$w(z)=u(2(z-c_n))$ in $\{z:|z-c_n|\le e^{n-1}\}\setminus I_n$,
and deduce
$$
\omega(i,I_n,\C\setminus E)\ge c\cdot e^{(\e-1)n}/n, \qquad n\ge 1. 
$$        
\end{proof}

\begin{rmk}
To estimate harmonic measure $\omega_E$ from above, we may use
Theorem~\ref{te2} or Theorems~D and E. In particular, in the
conditions of Lemma~{\rm\ref{le6}}, for some positive $c$ we have
$$
\omega(i,I_n,\C\setminus E)\le \exp(-c \sqrt n),\qquad n\ge 1.
$$
\end{rmk}

\noindent
{\bf Acknowledgment.} The authors thank 
Alexander~Fryntov for very useful discussions.

\textrm{\newline
\noindent
A.~Borichev
\newline
\noindent
Laboratoire de Math\'ematiques Pures de Bordeaux
\newline
\noindent
UPRESA 5467 CNRS, Universit\'e Bordeaux I
\newline
\noindent
351, cours de la Lib\'eration, 33405 Talence 
\newline
\noindent
FRANCE
\newline
\noindent
borichev@math.u-bordeaux.fr
\newline
\noindent
\newline
\noindent
M.~Sodin
\newline
\noindent
School of Mathematical Sciences 
\newline
\noindent
Tel-Aviv University
\newline
\noindent
Ramat-Aviv, 69978
\newline
\noindent
ISRAEL
\newline
\noindent
sodin@math.tau.ac.il}
\bigskip


\bigskip

\begin{thebibliography}{2}

\bibitem{Ak}
N.~I.~Akhiezer,
\newblock{\it On the weighted approximation of continuous functions
by polynomials on the entire real axis,}
\newblock{Uspekhi Mat.\ Nauk {\bf 11} (1956), 3--43; English
translation in American Mathematical Society 
Translations (ser.~2) {\bf 22} (1962), 95--137.}

\bibitem{AL}
N.~I.~Akhiezer, B.~Ja.~Levin,
\newblock{\it Generalization of S.~N.~Bernstein's 
inequality for derivatives of entire functions,} 
\newblock{Issledovaniya po sovremennym problemam teorii 
funkcii kompleksnogo peremennogo, Gosudarstv.\ Izdat.\ Fiz.-Mat.\ 
Lit., Moscow 1960, 111--165. French translation in}
\newblock{Fonctions d'une variable complexe. Probl\`emes 
contemporains (ed. A. Marcouchevitch), Gauthier--Villars, Paris, 
1962, 109--161.}

\bibitem{BA} M.~Benedicks,
\newblock{\it Positive harmonic functions vanishing on the boundary of
certain domains in $\R^n$,}
\newblock{Arkiv f\"or Matematik {\bf 18} (1980), 53--72.}

\bibitem{BP} M.~Benedicks,
\newblock{\it Weighted polynomial approximation on subsets of the real line,}
\newblock{Preprint 1981:11, Uppsala Univeristy, Mathematical Department
(1981), 12 pp.}

\bibitem{BM} A.~Beurling, P.~Malliavin,
\newblock{\it On Fourier transforms of measures with compact support,}
\newblock{Acta Mathematica {\bf 107} (1962), 291--309.}

\bibitem{deB1} L. de Branges,
\newblock{\it The Bernstein problem,}
\newblock{Proceedings of the American Mathematical Society 
{\bf 10} (1959), 825--832.}

\bibitem{deB2} L. de~Branges,
\newblock{\it Hilbert spaces of entire functions,}
\newblock{Prentice-Hall, 1968.}

\bibitem{GK}
I.~C.~Gohberg, M.~G.~Krein, 
\newblock{\it Introduction to the theory of linear
non-selfadjoint operators,} 
\newblock{Translations of Mathematical Monographs {\bf 18},
American Mathematical Society, 1969.}

\bibitem{GO}
A.~A.~Goldberg, I.~V.~Ostrovski\u\i, 
\newblock{\it The distribution of values of meromorphic functions,}
\newblock{Nauka, Moscow, 1970 (in Russian).}

\bibitem{Ha}
G.~H.~Hardy, 
\newblock{\it On the function $P_\rho(x)$,}
\newblock{Quart.\ Journ.\ Math.\ Ser.\ 2 {\bf 37} (1906), 
146--172; Complete works, vol.\ 4.}

\bibitem{HK}
W.~K.~Hayman, P.~B.~Kennedy, 
\newblock{\it Subharmonic functions,} 
\newblock{Vol. I, Academic Press, London et al, 1976.}

\bibitem{Kar}
P.~P.~Kargaev,
\newblock{\it Existence of the Phragm\'en--Lindel\"of function
and some conditions for quasianalyticity,}
\newblock{Zap.\ Nauchn.\ Sem.\ LOMI {\bf 126} (1983), 97--108; English
translation in J. Soviet Math.\ {\bf 27} (1984), 2486--2495.}

\bibitem{Koo}
P.~Koosis,
\newblock{\it The Logarithmic Integral,}
\newblock{Vol.\ I, Cambridge University Press, Cambridge, 1988.}

\bibitem{Kre1}
M.~G.~Krein,
\newblock{\it A contribution to the theory
of entire functions of exponential type,}
\newblock{Izv.\ Akad.\ Nauk SSSR {\bf 11} (1947), 309--326
(in Russian).}

\bibitem{Kre2}
M.~G.~Krein,
\newblock{\it Fundamental aspects of the representation theory  of
Hermitian operators with deficiency index $(m,m)$,} 
\newblock{Ukr.\ Math.\ Zh.\ {\bf 1} (1949), no.~2, 3--66; English 
translation in American Mathematical Society 
Translations (ser.~2), {\bf 97} (1970), 75--143.}

\bibitem{Kre3}
M.~G.~Krein,
\newblock{\it On the indeterminate case of the Sturm--Liouville 
problem in the interval $(0,\infty)$,}
\newblock{Izv.\ Akad.\ Nauk SSSR {\bf 16} (1952), 293--324
(in Russian).}

\bibitem{Lev2}
B.~Ya.~Levin,
\newblock{\it Density of systems of functions, quasianalyticity
and subharmonic majorants,}
\newblock{Zap.\ Nauchn.\ Seminarov LOMI {\bf 170} (1989), 102--156;
English translation in J. Soviet Math.\ {\bf 63} (1993), 171--201.}

\bibitem{LTF}
B.~Ya.~Levin,
\newblock{\it Majorants in classes of subharmonic functions, 
{\rm II}, The relation between majorants and conformal mapping,
{\rm III}, The classification of the closed sets on $\R$ and 
the representation of the majorants,}
\newblock{Teor.\ Funktsii Funktsional.\ Anal.\ i Prilozhen.\  
{\bf 52} (1989), 3--33;
English translation in J. Soviet Math.\ {\bf 52} (1990), 3351--3372.}

\bibitem{Lev1}
B.~Ya.~Levin,
\newblock{\it Lectures on entire functions,} 
\newblock{Translations of Mathematical Monographs {\bf 150},
American Mathematical Society, 1996.}

\bibitem{Me}
S.~N.~Mergelyan,
\newblock{\it Weighted approximation by polynomials,}
\newblock{Uspekhi Mat.\ Nauk {\bf 11} (1956), 107--152; English
translation in American Mathematical Society Translations (ser.~2) 
{\bf 10} (1958), 59--106.}

\bibitem{Nev}
R.~Nevanlinna, 
\newblock{\it Analytic functions,}
\newblock{Springer-Verlag, New York--Berlin, 1970.}

\bibitem{Par}
M.~Parreau, 
\newblock{\it Sur les moyennes des fonctions harmoniques et 
analytiques et la classification des surfaces de Riemann,} 
\newblock{Ann.\ Inst.\ Fourier Grenoble {\bf 3} (1951), 103--197.}

\bibitem{Sch}
A.~C.~Schaeffer, 
\newblock{\it Entire functions and trigonometric polynomials,} 
\newblock{Duke Math.\ Journ.\ {\bf 20} (1953), 77--88.}

\bibitem{SY}
M.~Sodin and P.~Yuditski\u\i,
\newblock{\it Another approach to de~Branges' theorem on weighted polynomial
approximation,} 
\newblock{in: Proceedings of the Ashkelon Workshop on Complex Function Theory 
(May 1996), L.~Zalcman, ed., Israel Mathematical Conferences 
Proceedings {\bf 11}, American Mathematical Society, 
Providence RI, 1997, pp. 221-227.} 

\end{thebibliography}
\end{document}